\documentclass[12pt, leqno]{amsart}
\usepackage{amsfonts,amssymb}
\usepackage{amsmath}
\let\geq\geqslant
\let\leq\leqslant

\theoremstyle{plain}
\newtheorem{theorem}{Theorem}[section]
\newtheorem{proposition}[theorem]{Proposition}
\newtheorem{lemma}[theorem]{Lemma}
\newtheorem{corollary}[theorem]{Corollary}

\theoremstyle{definition}
\newtheorem{definition}[theorem]{Definition}
\newtheorem{example}[theorem]{Example}
\newtheorem{remark}[theorem]{Remark}

\newcommand{\N}{\mathbb{N}}                   % natural numbers
\newcommand{\R}{\mathbb{R}}                   % real numbers
\newcommand{\C}{\mathcal{C}}
\newcommand{\sX}{\mathfrak{X}}                % stratification1
\newcommand{\W}{\mathcal{W}}                  % Whiteny condition
\newcommand{\Q}{\mathcal{Q}}
\newcommand{\WL}{\mathcal{WL}}

\begin{document}
\baselineskip=16pt
\title[Invariance of regularity conditions under definable, locally ... ]{Invariance of regularity conditions under definable, locally Lipschitz, weakly bi-Lipschitz mappings}
\author{Ma\l gorzata Czapla}
\address{Uniwersytet Jagiello\'{n}ski, Instytut Matematyki, ul. prof. Stanis\l awa \L ojasiewicza 6, 30-348 Krak\'ow, Poland}
\email{Malgorzata.Czapla@im.uj.edu.pl}

\begin{abstract}
\baselineskip=16pt

In this paper we describe the notion of a weak lipschitzianity of
a mapping on a $C^{q}$\ stratification. We also distinguish a class of regularity
conditions that are in some sense invariant under definable, locally Lipschitz and weakly
bi-Lipschitz homeomorphisms. This class includes the Whitney
(B) condition and the Verdier condition.
\end{abstract}
\maketitle

\section*{Introduction}
\label{sec: introduction}

The first goal of this paper is to study the notion of a weakly Lipschitz mapping on a fixed $C^{q}$\ stratification, which generalizes the notion
of a Lipschitz function. Section \ref{sec: prelimin} consists of the basic definitions and notation, while in Section \ref{sec: Weakly -lipschitz }
we introduce the main idea together with its geometrical interpretation and discuss its fundamental properties.
Weakly Lipschitz mappings were used earlier under the name \textit{fonctions rugueuses}\ by J.-L. Verdier \cite{Ver}, who considered them from a different
point of view compared to the present paper.

The second goal of this paper is to distinguish a special class of regularity conditions (Section \ref{sec: distinguished class of reg conditions});
namely, the conditions which are definable, generic, $C^{q}$\ invariant
and having the property of lifting with respect to locally Lipschitz mappings and the property of projection with respect to weakly Lipschitz mappings.
In the definable case these properties make the regularity conditions in some sense invariant with respect to definable, locally Lipschitz, weakly bi-Lipschitz
homeomorphisms (Theorem \ref{trm: Q invariance}).

In Sections \ref{sec: Whitney cond as Q conditions} and \ref{sec: Verider cond as Q cond} we prove that the Whitney (B) condition and the Verdier condition
belong to the distinguished class.

\section{Preliminaries}
\label{sec: prelimin}

We denote by $|\cdot|$\ the euclidean norm of $\R^{n}$, $S^{n-1}=\{x\in\R^{n}: |x|=1\}$. In the whole paper $q$\ denotes the class of smoothness of a mapping, so
$q\in\N$\ or $q\in\{\infty, \omega\}$, unless otherwise indicated. Also we denote by $\mathbb{G}_{k,n}$\ the Grassmann manifold of $k$\ dimensional vector subspaces of $\R^{n}$. Then $\mathbb{P}_{n-1}=\mathbb{G}_{1,n}$\ is the real projective space of dimension $n-1$.

\begin{definition}
\label{def: d(v,W)}
Let $v\in S^{n-1}$\ and $W$\ be a nonzero linear subspace of $\R^{n}$. We put
$$d(v, W)=\inf\{sin(v,w): w\in W\cap S^{n-1}\},$$
where $sin(v,w)$\ denotes the sine of the angle between the vectors $v$\ and $w$. We also put $d(u,W)=1$\ if $W=\{0\}$.
\end{definition}

\begin{definition}
\label{Funkcja d}
For any $P\in\mathbb{G}_{k,n}$\ and $Q\in\mathbb{G}_{l,n}$, we put
$$d(P,Q)=\sup\{d(\lambda; Q): \lambda\in P\cap S^{n-1}\}, $$
when $k>0$, and $d(P,Q)=0$, when $k=0$.
\end{definition}

Now we list some elementary properties of the function $d$, leaving the proof to the reader.

\begin{proposition}\label{prop: wlasnosci d}\ \\
a) $0\leq d(P,Q)\leq 1$. \\
b) $d(P,Q)=0 \Longleftrightarrow P\subset Q$. \\
c) $d(P,Q)=1 \Longleftrightarrow P\cap Q^{\perp}\neq \{0\}$\ (or $d(P,Q)<1 \Longleftrightarrow P\cap Q^{\perp}= \{0\}$\ ). \\
d) $d(P,R)\leq d(P,Q)+d(Q,R)$.\\
e) $d(P,Q)=d(Q,P)$\ if $\dim P=\dim Q$.\\
f) $d$\ is a metric on every $\mathbb{G}_{k,n}$.\\
g) $d(\R v, Q)=|v-\pi_{Q}(v)|=dist (v,Q)=|\pi_{Q^{\perp}}(v)|=\sin \left(v, \frac{\pi_{Q}(v)}{|\pi_{Q}(v)|}\right)=d\left(\R v, \R \frac{\pi_{Q}(v)}{|\pi_{Q}(v)|}\right)$,
where $\pi_{Q}$\ denotes the orthogonal projection onto $Q$, $v$\ is a unit vector, not orthogonal to $Q$\ and $dist(v,Q)=\inf \{|v-w|: w\in Q\}$.\\
h) Consider the following metric on $\mathbb{P}_{n-1}$:
\begin{center}
 $\widetilde{d}(\R v, \R w)=\min \{|u-w|, |u+w|\}$\hskip 3 mm for \hskip 2 mm $u,w\in S^{n-1}$.\\
 \end{center}
Then we have the following inequalities
$$\frac{1}{\sqrt{2}}\text{ }\widetilde{d}(\R v, \R w)\leq d(\R v, \R w)\leq \widetilde{d}(\R v, \R w).$$
i) If $P\subset P'$, then $d(P,Q)\leq d(P',Q).$ \\
j) If $Q'\subset Q$, then $d(P,Q)\leq d(P,Q').$ \\
\end{proposition}

To transform $d$\ into a metric $D$\ in $\mathbb{G}_{n}=\overset{n-1}{\underset{k=1}{\bigcup}}\mathbb{G}_{k,n}$\ (disjoint union), we put
$$D(P,Q)=\max \{d(P,Q), d(Q,P)\}.$$
Then $\mathbb{G}_{k,n}$\ are open-closed components in $\mathbb{G}_{n}$\ and $D(\mathbb{G}_{k,n}, \mathbb{G}_{l,n})=1$\ for $k\neq l$.
It is easy to check that
$$|d(P_{1},Q_{1})-d(P_{2},Q_{2})|\leq D(P_{1},P_{2})+D(Q_{1},Q_{2}),$$
hence the function $d$\ is continuous.

We will need another function which characterizes the mutual position of two linear subspaces $V$\ and $W$\ of $\R^{n}$.

\begin{definition}\label{def: d inf(V,W)}
 $$\delta(V,W)= \inf\{d(v, W): v\in V\cap S^{n-1}\}$$
if $V\neq \{0\}$, and $\delta(V,W)=1$\ if $V=\{0\}$.
\end{definition}

The reader will easily check the following properties

\begin{proposition}\label{prop: properties of d inf}\ \\
$i)$\ $\delta(V,W)=0 \Longleftrightarrow V\cap W \neq \{0\}$.\\
$ii)$\ $\delta(V,W)>0\Longleftrightarrow V\cap W = \{0\}$.\\
$iii)$\ $\delta(V,W)=1\Longleftrightarrow V\bot W$.\\
$iv)$\ $\delta(V,W)\leq d(V,W)\leq D(V,W)$\ if $V\neq \{0\}\neq W$.\\
$v)$\ $\delta$\ is continuous.\\
\end{proposition}

\begin{proposition}\label{prop: d inf, transv,separated Txf} Let $\Lambda\subset\R^{n}$\ be a $C^{q}$\ submanifold,
$f: \Lambda\longrightarrow\R^{m}$\ be a $C^{q}$\ mapping. Assume that for each $x\in graph\text{ }f|_{\Lambda}$
$$\delta(T_{x}graph\text{ }f|_{\Lambda}, \{0\}\times\R^{m})\geq\alpha>0,$$
where $\alpha$\ is a positive constant.
Then

$i)$\ for any point $x_{0}\in \overline{graph f|_{\Lambda}}$\ and for any sequence
$\{x_{\nu}\}_{\nu\in\N}\subset graph f|_{\Lambda}$\ converging to $x_{0}$\ such that the sequence $\left\{T_{x_{\nu}}graph f|_{\Lambda}\right\}_{\nu\in\N}$\
is convergent,
$$\delta\left(\underset{\nu\to +\infty}{\lim} T_{x_{\nu}}graph\text{ }f|_{\Lambda}, \{0\}\times\R^{m}\right)\geq\alpha>0.$$
\vskip 1 mm
$ii)$\ for any $C^{q}$\ submanifold $M\subset \Lambda$, $M\times\R^{m}$\ is transversal to $graph f|_{\Lambda}$\ in $\Lambda\times\R^{n}$.
\end{proposition}

\begin{proof}
Observe that $i)$\ follows from the continuity of $\delta$.
\vskip 1 mm
$ii)$\ By Proposition \ref{prop: properties of d inf} $ii$)\
$$T_{x}graph\text{ }f|_{\Lambda}\cap \left(\{0\}\times\R^{m}\right)=\{0\},$$
hence it is enough to observe that $\dim T_{x}graph f|_{\Lambda}=\dim\Lambda $.

\end{proof}

\begin{proposition}\label{prop: lipschitz mapping is separated from the vertical space} Let $\Lambda\subset\R^{n}$\ be a $C^{q}$\ submanifold and let
$f: \Lambda\longrightarrow\R^{m}$\ be a Lipschitz $C^{q}$\ mapping. Then there exists a positive constant $\alpha$\ such that
$$\delta(T_{x}graph\text{ }f, \{0\}\times\R^{m})\geq\alpha>0,$$
for each $x\in graph f$.
\end{proposition}
\begin{proof} \hskip 1 mm
There exists $L>0$\ such that $||d_{y}f||\leq L$, for each $y\in\Lambda$.
Now let $v\in T_{x}graph f$, $|v|=1$\ for some point $x=(y, f(y))$, $y\in\Lambda$. Then there exists a vector
$\widetilde{v}\in T_{y}\Lambda$\ such that
$$v=(\widetilde{v},d_{y}f(\widetilde{v})).$$
Then we have
$$1=|v|=\sqrt{|\widetilde{v}|^{2}+|d_{y}f(\widetilde{v})|^{2}}\leq \sqrt{1+L^{2}}\cdot |\widetilde{v}|.$$
Therefore
$$ |\widetilde{v}|\geq \frac{1}{\sqrt{1+L^{2}}}>0.$$
On the other hand
$$d(v, \{0\}\times \R^{n})=|\widetilde{v}|.$$
\end{proof}

Now we recall briefly the notion of a $C^{q}$\ stratification.

\begin{definition} Let $A$\ be a subset of $\R^{n}$.
\textit{A} $C^{q}$\ \textit{stratification} \textit{of} the set $A$\ is a (locally) finite family
$\sX_{A}$\ of connected $C^{q}$\ submanifolds of $\R^{n}$\ (called \textit{strata})\
such that \vskip 1 mm
\begin{tabular}{l}
$1)$\ \ $A=\bigcup\sX_{A}$\ ;\\
$2)$\ \ if\ $\Gamma_{1}, \Gamma_{2}\in\sX_{A}$, $\Gamma_{1}\neq\Gamma_{2}$\ then $\Gamma_{1}\cap\Gamma_{2}=\emptyset$; \\
$3)$\ \ for each $\Gamma\in\sX_{A}$\ the set
$(\overline{\Gamma}\setminus\Gamma)\cap A$\ is a union of
 some strata from \\
\hskip 6 mm the family $\sX_{A}$\ of dimension $<\dim\Gamma$.\\
\end{tabular}
\vskip 1 mm We say that the stratification $\sX_{A}$\ is
\textit{compatible with a family of sets} \linebreak $B_{i}\subset A$, $i\in I$\ if every set $B_{i}$\ is
a union of some strata of $\sX_{A}$.
\end{definition}

Actually, we will be interested only in finite stratifications.

\begin{definition}
Let $A\subset \R^{n}$\ and let $f: A\longrightarrow \R^{m}$\ be a continuous mapping, $\sX_{A}$\ be a $C^{q}$\ stratification of the
set $A$\ such that $f|_{\Gamma}$\ is of class $C^{q}$\ for all $\Gamma\in\sX_{A}$.
Then by the \emph{induced} $C^{q}$\ \emph{stratification} of the $graph f$, we will mean the following:
$$\sX_{graph f}(\sX_{A})=\{graph f|_{\Gamma}:\quad \Gamma\in\sX_{A}\}.$$
\end{definition}

A natural setting for our results is the theory of o-minimal structures (or more generally geometric categories), as presented in \cite{D} (or \cite{DM}).
In the whole paper the adjective \textit{definable} (i.e. definable subset, definable mapping) will refer to any fixed o-minimal structure on
the ordered field of real numbers $\mathbb{R}$.

\section{Weakly Lipschitz mappings}
\label{sec: Weakly -lipschitz }

In this section we describe the idea of the weak lipschitzianity of a mapping
and list its important properties.

\begin{definition}
\label{def:weakly_lipschitz-map} Let $A$\ be a subset of $\R^{n}$\ and let $\sX_{A}$\ be a finite $C^{q}$\
stratification of the set $A$. Consider a mapping $f: A\longrightarrow \R^{m}$.
\vskip 1 mm
We say that $f$\ is \textit{weakly Lipschitz of class }$C^{q}$\ \textit{on the stratification}\ $\sX_{A}$,
\ if for each stratum $\Gamma\in\sX_{A}$\ the restriction $f|_{\Gamma}$\ is of class $C^{q}$\ and the pair $(f, \sX_{A})$
\ satisfies one of the following equivalent conditions:
\vskip 1 mm
a) Whenever $\Lambda, \Gamma\in\sX_{A}$, $\Gamma\subset\overline{\Lambda}\setminus\Lambda$, $a\in\Gamma$\ and
$\{a_{\nu}\}_{\nu\in\N}$, $\{b_\nu\}_{\nu\in\N}$\ are arbitrary sequences such that $a_{\nu}\in \Gamma, b_{\nu}\in\Lambda$, for
$\nu\in\N$, then
$$a_{\nu}, b_{\nu}\longrightarrow a \quad(\nu\rightarrow +\infty)\quad \Longrightarrow
\quad \limsup_{\nu\rightarrow +\infty}\frac{|f(a_{\nu})-f(b_{\nu})|}{|a_{\nu}-b_{\nu}|}< +\infty. $$

b) For any stratum $\Gamma\in\sX_{A}$\ and any point $a\in\Gamma$\
there exists a neighbourhood $U_{a}$\ of $a$\ such that the
mapping
$$\psi: (\Gamma\cap U_{a})\times ((A\setminus \Gamma)\cap U_{a})\ni(x,y)\longmapsto \frac{|f(x)-f(y)|}{|x-y|}\in\R$$
is bounded.
\vskip 1 mm

c) For any strata $\Lambda, \Gamma\in\sX_{A}$, $\Gamma\subset\overline{\Lambda}\setminus\Lambda$\ and for any $a\in\Gamma$
\ there exists a neighbourhood $U_{a}$\ of $a$\ such that the mapping
$$\psi: (\Gamma\cap U_{a})\times (\Lambda\cap U_{a})\ni (x,y)\longmapsto \frac{|f(x)-f(y)|}{|x-y|}\in\R $$
is bounded.
\vskip 1 mm
d) Whenever $\Lambda, \Gamma\in\sX_{graph f}(\sX_{A})$, $\Gamma\subset\overline{\Lambda}\setminus\Lambda$,
$x\in\Gamma$\ and $\{x_{\nu}\}_{\nu\in\N}\subset \Gamma$, $\{y_\nu\}_{\nu\in\N}\subset\in\Lambda$\ are arbitrary sequences convergent to $x$, then
$$d\left(\underset{\nu\rightarrow +\infty}{\lim} \R\big(x_{\nu}-y_{\nu}),\text{  }\{0\}\times \R^{m}\right)>0.$$
\vskip 1 mm
e) Whenever $\Lambda, \Gamma\in\sX_{graph f}(\sX_{A})$, $\Gamma\subset\overline{\Lambda}\setminus\Lambda$, $x\in\Gamma$\ if
$\{x_{\nu}\}_{\nu\in\N}\subset \Gamma$, $\{y_\nu\}_{\nu\in\N}\subset\Lambda$\ are arbitrary sequences convergent to $x$\
and there exists a limit
$$L=\underset{\nu\rightarrow +\infty}{\lim} \R(x_{\nu}-y_{\nu}),$$
then $L\cap (\{0\}\times\R^{m})=\{0\}$.
\end{definition}

\begin{proposition}
The conditions a)-e)\ from Definition \ref{def:weakly_lipschitz-map} are equivalent.
\end{proposition}

\begin{proof}
$a)\Rightarrow c)$. Suppose that $c)$\ is not satisfied. Then we find some strata $\Gamma, \Lambda\in\sX_{A}$,
$\Gamma\subset\overline{\Lambda}\setminus\Lambda$\ and a point $a\in\Gamma$,\ for which we can find a basis of neighbourhoods
$\{U_{n}\}_{n\in\N}$\ and two sequences
$\{a_{n}\}_{n\in\N}$, $\{b_{n}\}_{n\in\N}$\ such that $a_{n}\in \Gamma\cap U_{n}$, $b_{n}\in\Lambda\cap U_{n}$\ and
$$\frac{|f(a_{n})-f(b_{n})|}{|a_{n}-b_{n}|}> n,$$
a contradiction with the assumption.
\vskip 1 mm
$c)\Rightarrow a)$, $b)\Leftrightarrow c)$\ are trivial.
\vskip 1 mm
$a)\Leftrightarrow d)\Leftrightarrow e)$. Since
$\Gamma=graph f|_{\Gamma'}$, $\Lambda=graph f|_{\Lambda'}$,  $x=(a, f(a))$, $x_{\nu}=(a_{\nu},f(a_{\nu}))$\ and $y_{\nu}=(b_{\nu},f(b_{\nu}))$,
when $\Gamma', \Lambda'\in \sX_{A}$, $a, a_{\nu}\in\Gamma'$, $b_{\nu}\in\Lambda'$, $a_{\nu}\longrightarrow a$\ and $b_{\nu}\longrightarrow a$\
($\nu\longrightarrow +\infty$), it is enough to observe that
%\begin{center}\begin{tabular}{l l}
{\setlength\arraycolsep{2pt}
\begin{eqnarray*}
d\left(\R(x_{\nu}-y_{\nu}), \{0\}\times\R^{m}\right)&=&\frac{|a_{\nu}-b_{\nu}|}{|(a_{\nu}, f(a_{\nu}))-(b_{\nu},f(b_{\nu}))|}\\
&=&\frac{1}{\sqrt{1+{\left(\frac{|f(a_{\nu})-f(b_{\nu})|}{|a_{\nu}-b_{\nu}|}\right)}^{2}}}.\end{eqnarray*}}

%\end{tabular}\end{center}
\end{proof}

\begin{remark} If $f: A\longrightarrow \R^{m}$\ is weakly Lipschitz on a stratification $\sX_{A}$\ of the set $A$, then
$f$\ is continuous on $A$.
\end{remark}

Of course, the weak lipschitzianity is a generalization of the Lipschitz condition. Obviously, we have the following

\begin{proposition}\label{prop: Lip=>weak-Lip}
Let $A\subset\R^{n}$, $f: A\longrightarrow \R^{m}$\ be a locally Lipschitz mapping.
Assume that $A$\ admits a $C^{q}$ stratification $\sX_{A}$\ such that for all strata $\Gamma\in\sX_{A}$\ the map $f|_{\Gamma}$\
is of class $C^{q}$. Then $f$\ is weakly Lipschitz of class $C^{q}$\ on the stratification $\sX_{A}$.
\end{proposition}

By the $C^{q}$\ Cell Decomposition Theorem (see \cite{DM}), we have the following

\begin{corollary} Let $A\subset\R^{n}$, $f: A\longrightarrow\R^{m}$\ be a definable locally Lipschitz mapping.
There exists a definable $C^{q}$\ stratification $\sX_{A}$\ of $A$\ such that $f$\ is weakly Lipschitz of class $C^{q}$\ on the stratification
$\sX_{A}$.
\end{corollary}

The weak lipschitzianity is a much weaker property than the local Lipschitz condition, as it is shown
in the examples below.

\begin{example}
Let $A\subset\R^{2}$, $A=\Lambda\cup\Gamma_{1}\cup\Gamma_{2}$\ and $\sX_{A}=\{\Lambda,\Gamma_{1},\Gamma_{2}\}$, where
\begin{center}
\begin{tabular}{l}
$\Lambda=\{(x,y)\in\R^{2}: x\in(0,1),\text{  } \frac{1}{2}x^{2}< y <x^{2}\}$, \\
$\Gamma_{1}=\{(x,y)\in\R^{2}: x\in(0,1),\text{  } y=\frac{1}{2}x^{2} \}$, \\
$\Gamma_{2}=\{(0,0)\}.$ \\
\end{tabular}
\end{center}
\vskip 2 mm
Consider the mapping $$f: A\ni(x,y)\longrightarrow (x,\sqrt{y})\in \R^{2}.$$
Then $f$\ is not Lipschitz in any neighbourhood of the point $(0,0)$, because $\frac{\partial f}{\partial y}=\left(0, \frac{1}{2\sqrt{y}}\right)$\ is unbounded. However,
$f$\ is weakly Lipschitz, because it is locally Lipschitz on $A\setminus\{(0,0)\}$\ and
$$\frac{|f(x,y)|}{|(x,y)|}=\sqrt{\frac{x^{2}+y}{x^{2}+y^{2}}}\leq \sqrt{\frac{2x^{2}}{x^{2}}}\leq \sqrt{2}.$$
\end{example}

\begin{example} Let $\Lambda=\{(x,y)\in\R^{2}:\text{ } y>0\}$, $\Gamma=\{(x,y)\in\R^{2}:\text{\ } y=0\}$, $A=\Lambda\cup\Gamma$.
Consider the mapping $f: A\longrightarrow \R$\ defined by the following formula
$$f(x,y)=\begin{cases}\left(\frac{x^{2}}{y^{7}}-y^{2}\right)^{2} & 0<x^{2}<y^{9}, \\
0 & x^{2}\geq y^{9}\geq 0. \end{cases}$$
\vskip 1 mm
Then $f$\ is weakly lipschitzian of class $C^{1}$\ on $\{\Lambda, \Gamma\}$.
Indeed, $f$\ is $C^{1}$\ on $A\setminus\{(0,0)\}$, hence it is locally Lipschitz on $A\setminus\{(0,0)\}$.
Moreover, if $(x,y)\in\Lambda$, $(x',0)\in\Gamma$\ and $f(x,y)\neq 0$, then
$$\frac{|f(x,y)-f(x',0)|}{|(x,y)-(x',0)|}=\frac{\left(\frac{x^{2}}{y^{7}}-y^{2}\right)^{2}}{\sqrt{(x-x')^{2}+y^{2}}}\leq \frac{y^{4}}{y}=y^{3}.$$
Nevertheless, $f$\ is not a Lipschitz mapping in any neighbourhood of $(0,0)$, because
$$\left|\frac{\partial f}{\partial x}\left(\frac{y^{\frac{9}{2}}}{\sqrt{3}},y\right)\right|=
\frac{8}{3\sqrt{3}}\frac{1}{\sqrt{y}}\longrightarrow +\infty,$$
if $y\longrightarrow 0$.
\end{example}

The proofs of the following three propositions are straightforward.

\begin{proposition}
\label{prop: substrat-preserv-WL} Let $A\subset\R^{n}$, $f: A\longrightarrow \R^{n}$\ be weakly Lipschitz of class $C^{q}$\ ($q\geq 1$)\
on a $C^{q}$\ stratification $\sX_{A}$. Let $B\subset A$.
Then for any $C^{q}$\ stratification $\sX_{B}$\ of the set $B$, compatible with $\sX_{A}$, the mapping $f$\ is weakly Lipschitz of
class $C^{q}$\ on the stratification $\sX_{B}$.
\end{proposition}

\begin{proposition}\label{prop: comp-weak-lip-weak-lip} Let $f: A\longrightarrow \R^{p}$\ be a weakly Lipschitz $C^{q}$\ mapping ($q\geq 1$)
on a $C^{q}$\ stratification $\sX_{A}$\ of a set $A\subset\R^{n}$\ and let $g: B\longrightarrow\R^{r}$\ be a weakly Lipschitz $C^{q}$\ mapping on a
$C^{q}$\ stratification $\sX_{B}$\ of a set $B\subset\R^{p}$. Assume that the image under $f$\ of each stratum from $\sX_{A}$\ is contained in some stratum
from $\sX_{B}$\ (in particular, $f(A)\subset B$). Then $g\circ f: A\longrightarrow \R^{r}$\ is a weakly Lipschitz $C^{q}$\ mapping on $\sX_{A}$.
\end{proposition}

\begin{remark}\label{rem: stratified sets WL category} The last proposition allows to define a category of stratified sets as objects and weakly Lipschitz
$C^{q}$\ mappings ($q\geq 1$)\ as morphisms.
\end{remark}

\begin{proposition}\label{prop: zest-WL-jest-WL} Let $f: A\longrightarrow \R^{m}$\ and $g: A\longrightarrow\R^{p}$\ be two weakly Lipschitz $C^{q}$\ mappings
on a $C^{q}$\ stratification $\sX_{A}$\ of a set $A\subset\R^{n}$. Then the mapping
$$(f,g): A\ni x\longmapsto (f(x),g(x))\in \R^{m}\times\R^{p}$$
is weakly Lipschitz of class $C^{q}$\ as well.
\end{proposition}

\begin{definition} For a homeomorphic embedding $f: A\longrightarrow \R^{m}$\ and a $C^{q}$\ stratification $\sX_{A}$\ of $A$\ such that
for any $\Gamma\in\sX_{A}$\ the map $f|_{\Gamma}$\ is a $C^{q}$\ embedding, we have a natural $C^{q}$\ stratification of the image $f(A)$\
$$f\sX_{A}=\{f(\Gamma):\quad \Gamma\in\sX_{A}\}.$$
\end{definition}

This leads to the following definition of a weakly bi-Lipschitz homeomorphism:

\begin{definition} Let $A\subset \R^{n}$\ be a set, $f:A\longrightarrow \R^{m}$\ be a homeomorphic embedding. Let
$\sX_{A}$\ be a $C^{q}$\ stratification ($q\geq 1$)\ of the set $A$\ such that for all $\Gamma\in\sX_{A}$\ the mapping $f|_{\Gamma}$\ is a
$C^{q}$\ embedding.

We say that the mapping $f$\ is \textit{weakly} \textit{bi-Lipschitz of class} $C^{q}$\ \textit{ on the stratification }$\sX_{A}$,\ if $f$\
is weakly Lipschitz of class $C^{q}$\ on $\sX_{A}$\ and the inverse mapping $f^{-1}: f(A)\longrightarrow A\subset\R^{n}$\ is weakly Lipschitz of class $C^{q}$\ on the
stratification $f\sX_{A}$.
\end{definition}

In order to check that the inverse mapping is weakly Lipschitz on $f\sX_{A}$, we can use the following obvious

\begin{proposition}
\label{prop: weak-lip-condition-for-inverse-mapping}
Let $A\subset \R^{n}$, $f: A\longrightarrow \R^{m}$\ be a homeomorphic embedding. Let $\sX_{A}$\ be a
$C^{q}$\ stratification of the set $A$\ and assume that for each stratum $\Gamma\in\sX_{A}$\ the mapping $f|_{\Gamma}$\ is
a $C^{q}$\ embedding.

Then the mapping $f^{-1}: f(A)\longrightarrow A$\ is weakly Lipschitz of class $C^{q}$\ on the stratification $f\sX_{A}$,\
if and only if it satisfies the following condition
\vskip 2 mm
a') for any strata $\Gamma, \Lambda\in \sX_{A}$, $\Gamma\subset\overline{\Lambda}\setminus\Lambda$\ and for any point $a\in\Gamma$\
if $\{a_{\nu}\}_{\nu\in\N}$, $\{b_\nu\}_{\nu\in\N}$\ are arbitrary sequences such that
$a_{\nu}\in \Gamma, b_{\nu}\in\Lambda$\ for $\nu\in\N$, then
$$a_{\nu}, b_{\nu}\longrightarrow a \quad(\nu\rightarrow +\infty)\quad \Longrightarrow \quad
\liminf_{\nu\longrightarrow +\infty}\frac{|f(a_{\nu})-f(b_{\nu})|}{|a_{\nu}-b_{\nu}|}> 0. $$
\end{proposition}

\section{The $\WL$\ class of regularity conditions}
\label{sec: distinguished class of reg conditions}
\vskip 3 mm In this section we describe some class of regularity conditions and we
prove that in o-minimal geometry they are in some sense invariant under definable, locally Lipschitz, weakly bi-Lipschitz homeomorphisms.
As it is shown in the next sections, this class includes the Whitney
(B) condition and the Verdier condition.

From now on we fix on the ordered field $\mathbb{R}$\ an o-minimal structure\ $\mathcal{D}$\ admitting definable $C^{q}$\ Cell Decompositions ($q\geq 1$\ is also fixed).
In the whole paper definable means definable in $\mathcal{D}$.

\begin{theorem} Let $p\in\N$, $p\geq 1$. For any finite family of definable sets $A$,
$B_{1},...,B_{p}\subset A\subset\R^{n}$, there exists a finite
definable $C^{q}$\ stratification of $A$, compatible with $B_{i}$,
$i=1,2,...,p$.
\end{theorem}
\begin{proof} This easily follows from $C^{q}$\ Cell Decomposition Theorem (see \cite{DM}, 4.2).

\end{proof}

Let $\Q$\ be a condition on pairs $(\Lambda,\Gamma)$\ at points $x\in\Gamma$, where $\Lambda,\Gamma\subset\R^{n}$\ are $C^{q}$\ submanifolds,
$\Gamma\subset\overline{\Lambda}\setminus\Lambda$. Sometimes we will refer to $\Q$\ as a \textit{regularity condition}.

\begin{definition} We say that $\mathcal{Q}$\ is \textit{local}\ if for an open neighbourhood $U$\ of the point $x\in\Gamma$
\ the pair $(\Lambda, \Gamma)$\ satisfies the condition $\Q$\ at $x$\ if and only if the pair $(\Lambda\cap U,\Gamma\cap U)$\ satisfies the condition $\Q$\ at the point $x$.
\end{definition}

We will be considering only local conditions.
We adopt the following notation:
\begin{center}
\begin{tabular}{l}
$\W^{\mathcal{Q}}(\Lambda,\Gamma,x)\equiv$\ the condition $\Q$\ is satisfied for the pair $(\Lambda,\Gamma)$\ at $x\in\Gamma$.\\
$\W^{\Q}(\Lambda,\Gamma)\equiv$\ $\W^{\mathcal{Q}}(\Lambda,\Gamma,x)$, for each $x\in\Gamma$. \\
$\sim\W^{\Q}(\Lambda, \Gamma,x)\equiv$\ the negation of $\W^{Q}(\Lambda,\Gamma,x)$.\\
\end{tabular}
\end{center}

\begin{definition}\label{def: defin cond}
We say that $\Q$\ is \textit{definable} if for any definable
$C^{q}$\ submanifolds $\Gamma,\Lambda\subset\R^{n}$\ such that $\Gamma\subset\overline{\Lambda}\setminus\Lambda$, the set
$$\{x\in\Gamma:\quad \W^{\Q}(\Lambda,\Gamma,x)\} $$
is definable.
\end{definition}

\begin{definition}\label{def: generic cond}
Let $\mathcal{Q}$\ be a definable condition. We say that $\mathcal{Q}$\ is \textit{generic} if for any definable
$C^{q}$\ submanifolds $\Lambda,\Gamma\subset\R^{n}$, $\Gamma\subset\overline{\Lambda}\setminus\Lambda$, the set
$$\{x\in\Gamma:\quad \sim\W^{\Q}(\Lambda,\Gamma,x)\}$$
is nowhere dense in $\Gamma$.
\end{definition}

We are interested in $C^{q}$\ stratifications satisfying a condition $\Q$.

\begin{definition}\label{def: Q - stratification} Let $\sX_{A}$\ be a $\C^{q}$\
stratification of $A$. We say that $\sX_{A}$\ is a \textit{stratification with a}\
\textit{condition}\ $\Q$\ (or $\Q$\ \textit{- stratification}) if for each pair $\Lambda,\Gamma\subset\sX_{A}$
$$\Gamma\subset\overline{\Lambda}\setminus\Lambda\qquad \Longrightarrow\qquad \W^{\Q}(\Lambda,\Gamma).$$
\end{definition}

\begin{theorem}\label{trm: Q-stratification} (\L ojasiewicz-Stasica-Wachta) Let $\Q$\ be a definable, generic condition.
Then for any finite family of definable subsets $\{A_{j}\}_{j\in I}$\ of $\R^{n}$,\ there exists a finite definable $C^{q}$\ stratification $\sX_{\R^{n}}$\ of
$\R^{n}$\ with the condition $\Q$, compatible with $\{A_{j}\}_{j\in I}$.
\end{theorem}
\begin{proof}
A procedure was given in \cite{LSW}, Prop.2 in the subanalytic case. It suffices to observe that the same argument works in a general definable
case.
\end{proof}

\begin{corollary}\label{col: Q strat dim A }Let $\Q$\ be a definable, generic condition. Given a definable $C^{q}$\ stratification $\sX_{A}$
\ of a definable set $A\subset\R^{n}$,\ there exists a finite definable $C^{q}$\ stratification $\sX_{A}'$\ of the set $A$\ with the
condition $\Q$\ such that $\sX_{A}'$\ is compatible with $\sX_{A}$\ and moreover
$$\{\Gamma'\in\sX_{A}:\quad \dim\Gamma'=\dim A\}=\{\Gamma\in\sX_{A}:\quad \dim\Gamma=\dim A\}.$$
\end{corollary}
\begin{proof} Observe that using the downward induction from \cite{LSW} and correcting a definable $C^{q}$\ stratification to the one satisfying
the condition $\Q$, all we need to do is to substratify these strata of $\sX_{A}$\ that are of dimension $< \dim A$.
\end{proof}

\begin{definition}\label{def: invariant condition Q} We say that a condition $\Q$\ is $C^{q}$\ \textit{invariant}
(or \textit{invariant under} $C^{q}$\ \textit{diffeomorphisms}) if for any $C^{q}$\ diffeomorphism $\phi: U\longrightarrow W$\ of open subsets $U,W\subset \R^{n}$\
and any $C^{q}$\ submanifolds $\Lambda$, $\Gamma\subset U$\ such that $\Gamma\subset\overline{\Lambda}\setminus\Lambda$
and for any point $x\in\Gamma$
$$\W^{Q}(\Lambda,\Gamma,x)\quad \Longleftrightarrow \quad \W^{\Q}(\phi(\Lambda), \phi(\Gamma),\phi(x)).$$
\end{definition}

The class of conditions that we are to describe consists of the conditions that are definable, generic and invariant
 under definable $C^{q}$\ diffeomorphisms. Additionally, two more features are required.

\begin{definition}\label{def: property proj weakly lipsch}
We say that a condition $\Q$\ has \textit{a projection property with respect to weakly Lipschitz mappings of class }$C^{q}$\
if for any $C^{q}$\ mapping $f: A\longrightarrow \R^{m}$\ weakly Lipschitz on a $C^{q}$\ stratification $\sX_{A}$\ of a set $A\subset\R^{n}$, we have
$$\sX_{graph\text{ }f}\left(\sX_{A}\right)\text{\ is a\ }\Q\text{\ -\ stratification}\quad\Longrightarrow \quad \sX_{A}\text{\ is a\ }\Q\text{\ -\ stratification}.$$
Notice that this condition is equivalent to the following one:

For any subset $E\subset\R^{n}\times \R^{m}$\ and its $C^{q}$\ stratification $\sX_{E}$\ such that $\pi_{1}|_{E}: E\longrightarrow\R^{n}$
\ is a homeomorphic embedding and for each $\Gamma\in\sX_{E}$, $\pi_{1}|_{\Gamma}: \Gamma\longrightarrow \R^{n}$\ is a $C^{q}$\ embedding\ and
$\left(\pi_{2}|_{E}\right)\circ \left(\pi_{1}|_{E}\right)^{-1}$\ is weakly Lipschitz of class $C^{q}$\ on $\pi_{1}\sX_{E}$, we have
$$\sX_{E}\ \text{\ is a\ }\Q\text{\ -\ stratification}\quad\Longrightarrow \quad \pi_{1}\sX_{E}\ \text{\ is a\ }\Q\text{\ -\ stratification},$$
where $\pi_{1}: \R^{n}\times\R^{m}\longrightarrow\R^{n}$\ and $\pi_{2}:\R^{n}\times\R^{m}\longrightarrow \R^{m}$\ denote natural projections.
\end{definition}

\begin{remark} By Propositions \ref{prop: comp-weak-lip-weak-lip} and \ref{prop: zest-WL-jest-WL} the condition that the mapping
$\left(\pi_{2}|_{E}\right)\circ \left(\pi_{1}|_{E}\right)^{-1}$\
is weakly Lipschitz of class $C^{q}$\ on $\pi_{1}\sX_{E}$\ is equivalent to the condition that $\pi_{1}|_{E}$\ is weakly bi-Lipschitz of class $C^{q}$\ on the stratification
$\sX_{E}$.
\end{remark}

A proof of the following proposition is trivial.

\begin{proposition}\label{prop: property proj weakly lipsch GENERAL VERSION}
Let $\Q$\ be a condition, having the projection property with respect to weakly Lipschitz mappings of class $C^{q}$.
Let $f: A\longrightarrow \R^{m}$\ be weakly Lipschitz of class $C^{q}$\ on a stratification $\sX_{A}$\ of a set $A\subset\R^{n}$.

Then for any $C^{q}$\ submanifolds $\Gamma,\Lambda\subset graph f$\ such that $\Gamma\subset\overline{\Lambda}\setminus\Lambda$, $\dim \Gamma<\dim \Lambda$\
\footnote{The inequality $\dim \Gamma<\dim\Lambda$\ is required for $\{\Lambda,\Gamma\}$\ to be a stratification of $\Lambda\cup\Gamma$. In the definable case
the inequality follows from $\Gamma\subset\overline{\Lambda}\setminus\Lambda$.}
 and
$\{\Lambda, \Gamma\}$\ are compatible with the stratification $\sX_{graph f}(\sX_{A})$
$$\W^{\Q}(\Lambda,\Gamma)\quad\Longrightarrow \quad \W^{\Q}(\pi_{1}(\Lambda),\pi_{1}(\Gamma)).$$
\vskip 1 mm
\end{proposition}

\begin{definition}\label{def: property lifting by Lipschitz}
We say that a condition $\Q$\ has \textit{a lifting property with respect to locally Lipschitz mappings of class }$C^{q}$\ if
for any two $C^{q}$\ submanifolds $\Lambda,\Gamma\subset\R^{n}$\ such that $\Gamma\subset\overline{\Lambda}\setminus\Lambda$,
and for any locally Lipschitz mapping $f: \Lambda\cup\Gamma\longrightarrow \R^{m}$\ such that the restrictions
$f|_{\Lambda}$, $f|_{\Gamma}$\ are of class $C^{q}$\ and for any $C^{q}$\ submanifolds $M,N\subset\R^{n}$\ such that
$N\subset\overline{M}\setminus N$\ and $\{M,N\}$\ is
  compatible with $\{\Lambda,\Gamma\}$, we have
  $$\W^{\Q}(M,N), \W^{\Q}(graph f|_{\Lambda},graph f|_{\Gamma}) \quad \Longrightarrow \quad \W^{\Q}(graph f|_{M},graph f|_{N}).$$
Equivalently, $\Q$\ has the lifting property with respect to locally Lipschitz mappings of class $C^{q}$\ if
for any two $C^{q}$\ submanifolds $\Lambda,\Gamma\subset\R^{n}\times \R^{m}$\ such that $\Gamma\subset\overline{\Lambda}\setminus\Lambda$\ and
$\pi_{1}|_{\Lambda\cup\Gamma}$\ is a homeomorphic embedding, $\pi_{1}|_{\Lambda}$, $\pi_{1}|_{\Gamma}$\ are $C^{q}$\ embeddings\ and
the mapping $\pi_{2}|_{\Lambda\cup\Gamma}\circ \left(\pi_{1}|_{\Lambda\cup\Gamma}\right)^{-1}: \pi_{1}(\Lambda)\cup\pi_{1}(\Gamma)\longrightarrow \R^{m}$
\ is locally Lipschitz, the following holds true:
\vskip 2 mm

 for any $C^{q}$\ submanifolds $M,N\subset\R^{n}$\ such that $N\subset\overline{M}\setminus M$,
{\setlength\arraycolsep{2pt}
\begin{eqnarray*}
i)& \W^{\Q}(M,N), \W^{\Q}(\Lambda,\Gamma), M\subset\pi_{1}(\Lambda),  N\subset\pi_{1}(\Gamma)\Longrightarrow\qquad\qquad \qquad\quad \\
& \qquad\qquad \qquad\qquad\qquad\qquad \qquad \W^{\Q}\left((M\times \R^{m})\cap \Lambda, (N\times\R^{m})\cap \Gamma\right) \\
ii)& \W^{\Q}(M,N), M,N\subset \pi_{1}(\Lambda) \Longrightarrow  \W^{\Q}((M\times \R^{m})\cap \Lambda, (N\times\R^{m})\cap \Lambda) \\
iii)& \W^{\Q}(M,N), M,N\subset \pi_{1}(\Gamma) \Longrightarrow  \W^{\Q}((M\times \R^{m})\cap \Gamma, (N\times\R^{m})\cap \Gamma).
\end{eqnarray*}
}
\end{definition}

\begin{remark} If a condition $\Q$\ is $C^{q}$\ invariant, then the implications $ii)$\ and $iii)$\
are always satisfied.
\end{remark}

\begin{definition}\label{def: WL condition} We say that a condition $\Q$\ is a $\mathcal{WL}$\ \textit{condition}\
\textit{of class}\ $C^{q}$\ if it is\\
- definable ; \\
- generic ; \\
- invariant under definable $C^{q}$\ diffeomorphisms ;\\
- has the projection property with respect to weakly Lipschitz mappings of class $C^{q}$; \\
- has the lifting property with respect to locally Lipschitz mappings of \linebreak class $C^{q}$.
\end{definition}

Conditions of type $\mathcal{WL}$\ are invariant under definable locally Lipschitz, weakly bi-Lipschitz homeomorphisms in the following sense:

\begin{theorem}(Invariance of $\WL$\ conditions under definable, locally Lipschitz, weakly bi-Lipschitz homeomorphisms)
\label{trm: Q invariance}\hskip 1 mm
Let $\Q$\ be a $\mathcal{WL}$\ condition of class $C^{q}$. Let $A\subset \R^{n}$\ be a definable set and let $f:A\longrightarrow \R^{m}$\
be a definable homeomorphic embedding, weakly bi-Lipschitz of class $C^{q}$\ on a definable $C^{q}$\ stratification $\sX_{A}$\ of the set $A$. Assume
that for any two submanifolds $\Lambda$, $\Gamma\in \sX_{A}$\ such that
$\Gamma\subset\overline{\Lambda}\setminus\Lambda$,\ the mapping $f|_{\Lambda\cup\Gamma}$\ is locally Lipschitz.
\vskip 2 mm

Then there exists a definable $C^{q}$\ stratification $\sX_{A}'$\ of $A$, compatible with $\sX_{A}$\ such that
 $$\{\Gamma\in\sX_{A}:\quad \dim\Gamma=\dim A\}=\{\Gamma'\in\sX_{A}':\quad \dim\Gamma'=\dim A\}$$
and such that the condition $\Q$\ is invariant with respect to the pair $(f,\sX_{A}')$\ in the following sense
\vskip 2 mm
 for any definable $C^{q}$\ submanifolds $M,N\subset A$\ such that $N\subset\overline{M}\setminus M$\ and
$\{M,N\}$\ are compatible with the stratification $\sX_{A}'$\
$$\W^{\Q}(M,N)\quad\Longrightarrow\quad \W^{\Q}(f(M),f(N)).$$
\end{theorem}

\begin{proof}
Consider $graph f\subset \R^{n}\times\R^{m}$\ and its definable $C^{q}$\ stratification
$$\sX_{graph f}(\sX_{A})=\{graph f|_{\Gamma}:\text{ }\Gamma\in\sX_{A}\}.$$

By Corollary \ref{col: Q strat dim A } we find a definable $C^{q}$\ substratification $\sX^{\Q}_{graph f}(\sX_{A})$\
of the $graph f$\ with the condition $\Q$\ that is compatible with the family $\sX_{graph f}(\sX_{A})$\ and moreover
$$\{\Gamma'\in \sX^{\Q}_{graph f}(\sX_{A}):\dim\Gamma'=\dim A\}=\{\Gamma\in\sX_{graph f}(\sX_{A}):\dim\Gamma=\dim A\}.$$

Now we will show that
$$\sX'_{A}=\{\pi_{1}(\Lambda):\text{  }\Lambda\in\sX^{\Q}_{graph f}(\sX_{A})\}$$
is a required stratification.
\vskip 2 mm

Of course $\sX'_{A}$\ is really a definable $C^{q}$\ stratification of $A$, compatible with $\sX_{A}$\ and such that
$$\{\Gamma\in\sX_{A}: \text{  } \dim \Gamma=\dim A\}=\{\Lambda\in\sX'_{A}:\text{  }\dim\Lambda= \dim A\}.$$
By the projection property of the condition $\Q$\ with respect to weakly Lipschitz mappings, $\sX'_{A}$\ is a $\Q $\ - stratification.
\vskip 1 mm

Observe that for any two strata $\Lambda,\Gamma\in\sX_{A}'$, $\Gamma\subset\overline{\Lambda}\setminus\Lambda$,\ the mapping
$f|_{\Lambda\cup\Gamma}$\ is still locally Lipschitz, because of the compatibility of $\sX_{A}'$\ with $\sX_{A}$.
\vskip 1 mm

For the same reason $f\sX'_{A}$\ is compatible with $f\sX_{A}$.
Therefore, from Proposition \ref{prop: substrat-preserv-WL} we get that  $f$\ is still weakly bi-Lipschitz of class $C^{q}$\ on
$\sX'_{A}$.
\vskip 1 mm

Consider now two definable $C^{q}$\ submanifolds $M, N$\ in $\R^{n}$, $N\subset\overline{M}\setminus M$\ that are
compa\-ti\-ble with the stratification $\sX'_{A}$\ and $\W^{\Q}(M,N)$. Then two cases are possible:
\vskip 2 mm

I. There exists a submanifold $\Lambda\in \sX_{A}'$\ such that
$M,N\subset\Lambda$. As $f|_{\Lambda}$\ is a definable $C^{q}$\ embedding, then by the invariance of the condition $\Q$\
under definable $C^{q}$\ diffeomorphisms we get that $\W^{\Q}(f(M),f(N))$.
\vskip 2 mm

II. There exist two different $\Lambda,\Gamma\in\sX_{A}'$\ such that $N\subset\Gamma$, $M\subset\Lambda$. Then we have
$\Gamma\subset\overline{\Lambda}\setminus\Lambda$\ and observe that $f|_{\Lambda\cup\Gamma}$\ is still locally Lipschitz.
\vskip 2 mm
In this case, by the construction of the stratification $\sX^{\Q}_{graph\text{ }f}(\sX_{A})$\ we know that
$\W^{\Q}(graph f|_{\Lambda},graph f|_{\Gamma})$, so as the condition $\Q$\ has the lifting property under
locally Lipschitz mappings of class $C^{q}$, we get $$\W^{\Q}(graph f|_{M},graph f|_{N}).$$

On the other hand the mapping $f^{-1}$\ is weakly Lipschitz on the definable $C^{q}$\ stratification
$\{f(\Lambda),f(\Gamma)\}$\ and thus on $\{f(M), f(N)\}$. Observe that
$$\Phi\left(graph f|_{M}\right)= graph f^{-1}|_{f(M)},\qquad \Phi\left(graph f|_{N}\right)= graph f^{-1}|_{f(N)},$$
where $\Phi:\R^{n}\times \R^{m}\ni (x,y)\longmapsto (y,x)\in\R^{m}\times\R^{n}$. Therefore, by the projection property of the condition $\Q$\
with respect to weakly Lipschitz mappings we have
$$\W^{\Q}(\pi_{2}(graph f|_{M}), \pi_{2}(graph f|_{N})).$$
But $\pi_{2}(graph f|_{M})=f(M)$, $\pi_{2}(graph f|_{N}))=f(N)$,
which completes the proof.

\end{proof}

\begin{remark} As a corollary from the proof of Theorem \ref{trm: Q invariance} we get that the definable $C^{q}$\ stratifications
$\sX'_{A}$\ and $\sX'_{f(A)}=\{f(\Gamma):\text{\ } \Gamma\in \sX'_{A}\}$\ are the \linebreak $\Q$\ - stratifications.
\end{remark}

\section{The Whitney (B) condition as a $\WL$\ condition}
\label{sec: Whitney cond as Q conditions}

In this section we will prove that the Whitney (B) condition is of type $\WL$\ of class $C^{q}$, $q\geq 1$.

\begin{definition}\label{def: Whitney (B)} Let $N$, $M$\ be $C^{q}$\ submanifolds of $\R^{n}$\ such that
$N\subset \overline{M}\setminus M$\ and let $a\in N$.
\vskip 1 mm

We say that the pair of strata ($M$, $N$)\ satisfies \textit{the Whitney (B)}\ \textit{condition at the point} $a$\ (notation: $\W^{B}(M,N,a)$)\ if
for any sequences $\{a_{\nu}\}_{\nu\in \N}\subset N$, $\{b_{\nu}\}_{\nu\in\N}\subset M$\
both converging to the point $a$\ and such that the sequence of the secant lines $\{\R(a_{\nu}-b_{\nu})\}_{\nu\in\N}$
\ converges to a line $L\subset \R^{n}$\ in $\mathbb{P}_{n-1}$\ and the sequence of the tangent spaces $\{T_{b_{\nu}}M\}_{\nu\in\N}$\
converges to a subspace $T\subset \R^{n}$\ in $\mathbb{G}_{\dim M,n}$, always $L\subset T$.

When the pair of $C^{q}$\ submanifolds $(M,N)$\ satisfies (respectively, does not satisfy) the Whitney (B) condition\ at a point $a\in N$, we write $\W^{B}(M,N,a)$\
(respectively $\sim\W^{B}(M,N,a)$). If for any point $a\in N$\ we have $\W^{B}(M,N,a)$, we write $\W^{B}(M,N)$.
\end{definition}

\begin{theorem}\label{trm: (B) is definable} The Whitney (B) condition is definable and generic.
\end{theorem}
\begin{proof}
The proof in \cite{TL2} for the structure $(\R, +,\cdot, \exp, (r)_{r\in\R})$\ remains valid for arbitrary o-minimal structures on $(\R,<, +,\cdot)$.
See also \cite{DM}.

\end{proof}

\begin{definition}
A definable $C^{q}$\ stratification $\sX_{A}$ of a definable set $A\subset\R^{n}$\ is called \textit{a Whitney stratification} if for any pair of strata
$\Gamma, \Lambda\in\sX_{A}$\ such that $\Gamma\subset \overline{\Lambda}\setminus \Lambda$, the condition $\W^{B}(\Lambda,\Gamma)$\ is satisfied.
\end{definition}

\begin{remark}\label{trm: Whitney B is C1 invariant} The Whitney (B) condition is invariant under $C^{1}$\ diffeomorphisms (see \cite{Tro}).
\end{remark}

In order to show that the Whitney (B) condition has the projection property with respect to weakly Lipschitz mappings of class $C^{q}$, it suffices to prove the following
theorem.

\begin{theorem}\label{trm: Whitney (B) projective property}
Let $\Lambda, \Gamma\subset\R^{n}$\ be $C^{q}$\ submanifolds\ such that $\Gamma\subset\overline{\Lambda}\setminus\Lambda$\ and $\dim\Gamma<\dim\Lambda$.
Consider a mapping $f: \Lambda\cup\Gamma \longrightarrow \R^{m}$\ weakly Lipschitz
of class $C^{q}$\ on the stratification $\{\Lambda, \Gamma\}$.
Then
$$\W^{B}(graph f|_{\Lambda},graph f|_{\Gamma})\Longrightarrow \W^{B}(\Lambda,\Gamma).$$
\end{theorem}

\vskip 1 mm
\begin{proof} Let $a\in\Gamma$\ and $\{a_{k}\}_{k\in\N}\subset \Gamma$, $\{b_{k}\}_{k\in\N}\subset \Lambda$\ be sequences converging to $a$\ and such that
$$\R(a_{k}-b_{k})\longrightarrow L,\quad T_{b_{k}}\Lambda\longrightarrow T,\qquad \text{when}\quad k\longrightarrow +\infty$$
with some $L\in\mathbb{P}_{n-1}$, $T\in\mathbb{G}_{\dim \Lambda, n}$. Then
$$(a_{k},f(a_{k}))\longrightarrow (a,f(a)), \quad (b_{k},f(b_{k}))\longrightarrow (a,f(a)),\qquad k\longrightarrow +\infty.$$
Without loss of generality we may assume that for $k\longrightarrow+\infty$\ we have
$$\R\left((a_{k},f(a_{k}))-(b_{k},f(b_{k}))\right)\longrightarrow L',\quad T_{(b_{k}, f(b_{k}))}graph\text{ }f|_{\Lambda}\longrightarrow T'
$$
with some $L'\in\mathbb{P}_{n+m-1}$, $T'\in\mathbb{G}_{\dim \Lambda,n+m}$. Since $\W^{B}(graph f|_{\Lambda},graph f|_{\Gamma})$,
then $L'\subset T'$.
Because $f$\ is weakly Lipschitz on $\{\Lambda,\text{  }\Gamma\}$, thus
$$\pi_{1}(L')=L.$$
\vskip 1 mm
On the other hand, as $f|_{\Lambda}$\ is of class $C^{q}$, then for any $k\in\N$
$$\pi_{1}\left(T_{(b_{k}, f(b_{k}))}graph\text{ }f|_{\Lambda}\right)=T_{b_{k}}\Lambda,$$
hence, by the continuity of $\pi_{1}$\
we get the inclusion $\pi_{1}(T')\subset T$. Consequently,
$$L=\pi_{1}(L')\subset\pi_{1}(T')\subset T.$$
\end{proof}

Now we deal with the lifting property with respect to locally Lipschitz mappings of class $C^{q}$\ for the Whitney (B) condition.
It will easily follow from a more general transversal intersection theorem for the Whitney (B) condition.

\begin{theorem}\label{trm: Whitney transveral theorem}
Let $\Lambda_{i}, \Gamma_{i}\text{\ } (i=1,2)$\ be two pairs of $C^{q}$\ submanifolds in $\R^{n}$\ such that
$\Gamma_{i}\subset\overline{\Lambda_{i}}\setminus\Lambda_{i}$\ and $\W^{B}(\Lambda_{i},\Gamma_{i})$.
Assume that $\Lambda_{1}\cap\Lambda_{2}$\ and $\Gamma_{1}\cap\Gamma_{2}$\ are $C^{q}$\ submanifolds such that
$\Gamma_{1}\cap\Gamma_{2}\subset\overline{\Lambda_{1}\cap\Lambda_{2}}$\ and for each $x_{0}\in\Gamma_{1}\cap\Gamma_{2}$\ and any sequence
$\{y_{\nu}\}_{\nu\in\N}\subset\Lambda_{1}\cap\Lambda_{2}$\ converging to $x_{0}$, we have
$$T_{y_{\nu}}\Lambda_{i}\longrightarrow S_{i}\text{\ \ }(i=1,2)\quad \Longrightarrow \quad \dim S_{1}\cap S_{2}=\dim \Lambda_{1}\cap\Lambda_{2}.$$
Then $\W^{B}(\Lambda_{1}\cap\Lambda_{2}, \Gamma_{1}\cap\Gamma_{2})\ \footnotemark.$
\end{theorem}
\footnotetext{If for any $G\in\{\Lambda_{1},\Gamma_{1}\}$\ and $K\in\{\Lambda_{2},\Gamma_{2}\}$\ the submanifolds $G$\ and $K$\ are transversal in $\R^{n}$,
then Theorem \ref{trm: Whitney (B) lifting property} follows from Lemma 4.2.2 in \cite{Te}.}

\begin{proof} Let $x_{0}\in \Gamma_{1}\cap \Gamma_{2}$. Consider two sequences
\begin{center}$\{x_{n}\}_{n\in\N}\subset \Gamma_{1}\cap \Gamma_{2},\qquad \{y_{n}\}_{n\in\N}\subset \Lambda_{1}\cap \Lambda_{2}$,
\end{center}
such that $x_{n}, y_{n}\longrightarrow x_{0}$\ for $n\longrightarrow +\infty$. Assume that $\R (x_{n}-y_{n})\longrightarrow L$\ and that the
following sequences
\begin{center}
$\left\{T_{y_{n}}\left(\Lambda_{1}\cap \Lambda_{2}\right)\right\}_{n\in\N},\qquad \left\{T_{y_{n}}\Lambda_{1}\right\}_{n\in\N},
\qquad\left\{T_{y_{n}}\Lambda_{2}\right\}_{n\in\N}$\
\end{center}
are convergent. Let $$T= \lim_{y_{n}\to x_{0}}T_{y_{n}}(\Lambda_{1}\cap \Lambda_{2}).$$
By the assumptions
$$L\subset \lim_{y_{n}\to x_{0}}T_{y_{n}}\Lambda_{1}\cap \lim_{y_{n}\to x_{0}}T_{y_{n}}\Lambda_{2}=
 \lim_{y_{n}\to x_{0}}T_{y_{n}}(\Lambda_{1}\cap \Lambda_{2})=T.$$

\end{proof}

\begin{theorem}
\label{trm: Whitney (B) lifting property}
Let $\Lambda$, $\Gamma$\ be $C^{q}$\ submanifolds of $\R^{n}$\ such that $\Gamma\subset\overline{\Lambda}\setminus \Lambda$.
Let $f: \Lambda\cup\Gamma\longrightarrow \R^{m}$\ be a locally Lipschitz mapping such that $f|_{\Lambda}$, $f|_{\Gamma}$\ are of class $C^{q}$.
Let $M$, $N$\ be $C^{q}$\ submanifolds of\ $\R^{n}$\ such that $N\subset\overline{M}\setminus M$\ and $\{M, N\}$\ are compatible with $\{\Lambda,\Gamma\}$.
Then
$$\W^{B}(M,N),\quad \W^{B}(graph f|_{\Lambda},graph f|_{\Gamma})\quad \Longrightarrow \W^{B}(graph f|_{M}, graph f|_{N}).$$
\end{theorem}
\begin{proof}
If $M,N\subset \Lambda$\ or $M,N\subset\Gamma$, then $\W^{B}(graph f|_{M},graph f|_{N})$\ holds true, because the Whitney (B) condition
is $C^{q}$\ invariant. Now let $M\subset\Lambda$, $N\subset \Gamma$. It is enough to use Theorem \ref{trm: Whitney transveral theorem}, where
$$\Lambda_{1}=graph f|_{\Lambda}, \quad \Gamma_{1}=graph f|_{\Gamma},\quad \Lambda_{2}=M\times\R^{m},\quad \Gamma_{2}=N\times \R^{m}.$$
The last assumption of Theorem \ref{trm: Whitney transveral theorem} is fulfilled, because $f$\ is locally Lipschitz (use Prop.
\ref{prop: d inf, transv,separated Txf}
and \ref{prop: lipschitz mapping is separated from the vertical space}).

\end{proof}

\begin{corollary} The Whitney (B) condition is a $\WL$\ condition of class $C^{q}$.
\end{corollary}

\begin{corollary} Theorem \ref{trm: Q invariance} holds true for the Whitney (B) condition.
\end{corollary}

\begin{remark} A similar theorem holds true in the analytic-geometric category defined in \cite{DM}, as the Whitney stratification theorem holds true
in this category.
\end{remark}

\section{The Verdier condition as a $\WL$\ condition}
\label{sec: Verider cond as Q cond}

We start with some preparation.

\begin{lemma}\label{lem: Lipsch_rzut_funkcji_d}
Let $V$\ be a linear subspace of\ \ $\R^{n}$, $\R^{n}=V\oplus V^{\perp}$. Let $0<\alpha\leq 1$\ be a constant and consider a set
$$B_{\alpha}=\{u\in S^{n-1}:\quad d(\R u,V^{\perp})\geq \alpha\}.$$

Then there exists\ \ $C_{\alpha}>0$\ such that:

$i)$\ for any $u\in B_{\alpha}$, $w\in B_{\alpha}$\ we have
$$d(\R \pi_{V}(u),\R\pi_{V}(w))\leq C_{\alpha}\cdot d(\R u, \R w),$$

$ii)$\ for any two linear subspaces $L,K\subset\R^{n}$\ such that $L\cap S^{n-1}\subset B_{\alpha}$\ and $K\cap S^{n-1}\subset B_{\alpha}$,
$$d\left(\pi_{V}(L),\pi_{V}(K)\right)\leq C_{\alpha}\cdot d(L,K).$$

\end{lemma}
\begin{proof} $i)$. Let $u,w \in B_{\alpha}$. Then

\begin{tabular}{l}
$|\pi_{V}(u)|=|u-\pi_{V^{\perp}}(u)|= d(\R u,V^{\perp})\geq \alpha$, \\
$|\pi_{V}(w)|=|w-\pi_{V^{\perp}}(w)|= d(\R w,V^{\perp})\geq \alpha$. \\
\end{tabular}\\
Hence
\begin{eqnarray*}
%\begin{tabular}{l l}
&\left| \frac{\pi_{V}(u)}{\left|\pi_{V}(u)\right|}- \frac{\pi_{V}(w)}{\left|\pi_{V}(w)\right|}\right| \quad =
\left|\frac{\pi_{V}(u)\cdot\left|\pi_{V}(w)\right|-\pi_{V}(w)\cdot\left|\pi_{V}(u)\right|}{|\pi_{V}(u)|\cdot|\pi_{V}(w)|}\right| = \\
&= \left|\frac{\pi_{V}(u)\cdot\left|\pi_{V}(w)\right| - \pi_{V}(u)\cdot \left|\pi_{V}(u)\right|+
\pi_{V}(u)\cdot \left|\pi_{V}(u)\right|-\pi_{V}(w)\cdot\left|\pi_{V}(u)\right|}{|\pi_{V}(u)|\cdot|\pi_{V}(w)|}\right| \\
&\leq \frac{1}{\alpha^{2}}\cdot \big| \pi_{V}(u)\cdot\left(\left|\pi_{V}(w)\right|- \left|\pi_{V}(u)\right|\right)+
 \left|\pi_{V}(u)\right|\cdot \left(\pi_{V}(u)- \pi_{V}(w)\right)\big| \\
& \leq \frac{1}{\alpha^{2}}\cdot \big( |\pi_{V}(u)|\cdot\big|\left|\pi_{V}(u)\right|- \left|\pi_{V}(w)\right|\big|+
 \left|\pi_{V}(u)\right|\cdot \left|\pi_{V}(u)- \pi_{V}(w)\right|\big) \\
&\leq \frac{2}{\alpha^{2}}\big(|\pi_{V}(u)-\pi_{V}(w)|\big)=\frac{2}{\alpha^{2}}\big|\pi_{V}(u-w) \big|\leq\frac{2}{\alpha^{2}}|u-w|.
\end{eqnarray*}
Similarly,
$$\left| \frac{\pi_{V}(u)}{\left|\pi_{V}(u)\right|}+ \frac{\pi_{V}(w)}{\left|\pi_{V}(w)\right|}\right|=
 \left| \frac{\pi_{V}(u)}{\left|\pi_{V}(u)\right|}- \frac{\pi_{V}(-w)}{\left|\pi_{V}(-w)\right|}\right|\leq \frac{2}{\alpha^{2}}|u+w|.$$
Finally,
\begin{eqnarray*}d(\R\pi_{V}(u), \R\pi_{V}(w))\leq\widetilde{d}(\R\pi_{V}(u), \R\pi_{V}(w))&\leq \frac{2}{\alpha^{2}}\cdot \widetilde{d}(\R u,\R v) \\
&\leq\frac{2\sqrt{2}}{\alpha^{2}}d(\R u, \R v).
\end{eqnarray*}
The $ii)$\ is an easy corollary from $i)$.

\end{proof}

We will also need the following definition of the sine of the angle between two linear subspaces.

\begin{definition}\label{def: angle between spaces} Let $S,K\subset\R^{n}$\ be linear subspaces. Then we define
$$\lambda(S,K)=\begin{cases}\inf\{d(\R u, \R w): u\in S\cap S^{n-1},w\in K\cap S^{n-1}, u,w\perp S\cap K\}, \\
\hskip 6,7 cm \text{ for }S\not\subset K\text{ and } K\not\subset S, \\
0, \hskip 6,3 cm \text{ for }S\subset K \text{\ \ or\ \ }K\subset S.\end{cases}$$
\end{definition}

\begin{remark}\label{rem: wlasnosci lambda} Notice that if $S\not\subset K$\ and $K\not\subset S$, then
$$\lambda(S,K)=\delta\left(S\cap (S\cap K)^{\perp}, K\cap (S\cap K)^{\perp}\right).$$
\end{remark}

\begin{proposition}\label{prop: lambda pierwsza nierownosc} Let $k,l,n\in\N$\ and let $S\in\mathbb{G}_{s,n}$, $K\in\mathbb{G}_{k,n}$\ be such that $\lambda(S,K)>0$.
Then

$i)$\ for any $v\in\R^{n}$, $|v|=1$ we have
$$d(\R v, S\cap K)\leq \frac{1}{\lambda(S,K)}(d(\R v,S)+d(\R v,K)).$$

$ii)$\ If $R, L$\ are linear subspaces in $\R^{n}$, then
$$d(R\cap L,S\cap K)\leq \frac{1}{\lambda(S,K)}(d(R,S)+d(L,K)).$$
\end{proposition}
\begin{proof} $i)$. If $v\in S\cap K$, then the above inequality is satisfied as $$d(v,S\cap K)=d(v,S)=d(v,K)=0.$$

Assume now that $v\not\in S\cap K$, which means that $d(v,S\cap K)>0$. Then
$$\lambda(S,K)\leq d(\R (\pi_{K}(v)-\pi_{S\cap K}(v)), \R(\pi_{S}(v)-\pi_{S\cap K}(v)))$$
$$\hskip 1 cm \leq d(\R (\pi_{K}(v)-\pi_{S\cap K}(v)), \R(v-\pi_{S\cap K}(v)))+ $$
$$\hskip 5 cm d(\R (v-\pi_{S\cap K}(v)), \R(\pi_{S}(v)-\pi_{S\cap K}(v)))$$
$$=\frac{|v-\pi_{K}(v)|}{|v-\pi_{S\cap K}(v)|}+\frac{|v-\pi_{S}(v)|}{|v-\pi_{S\cap K}(v)|}=\frac{d(v,K)}{d(v,S\cap K)}+ \frac{d(v,S)}{d(v, S\cap K)},$$
as required.

The assertion $ii)$\ follows trivially from $i)$.

\end{proof}

\begin{proposition}
Let $p,k,s,n\in\N, p<\min\{k,s\}$.  Consider a set
$$\Sigma=\{(S,K)\in\mathbb{G}_{s,n}\times\mathbb{G}_{k,n}:\quad \dim (S\cap K)=p\}.$$
Then the mapping
$$\lambda: \Sigma \ni(S,K)\longmapsto \lambda(S,K)\in [0,1]$$
is continuous.
\end{proposition}
\begin{proof}
The continuity of $\lambda$\ follows easily from the fact that the mapping
$$\psi: \Sigma \ni (S,K)\longmapsto S\cap K\in\mathbb{G}_{p,n}$$
is continuous at any point $(S_{0}, K_{0})\in\Sigma$, because Proposition \ref{prop: lambda pierwsza nierownosc}$ii)$\
implies the inequality
$$d(S\cap K, S_{0}\cap K_{0})\leq \frac{1}{\lambda(S_{0},K_{0})}(d(S,S_{0})+d(K,K_{0}))$$
for any $(S,K)\in\Sigma$.

\end{proof}

\begin{proposition}\label{prop: inf lambda} Let $s,k,p,n\in \N$\ and $p<\min\{k,s\}$. Let $\widetilde{\Sigma}$\ be a compact subset of
the set
$$\Sigma=\{(S,K)\in\mathbb{G}_{s,n}\times\mathbb{G}_{k,n}:\quad \dim (S\cap K)=p\}.$$

Then
$$\inf\{\lambda(S,K):\quad (S,K)\in\widetilde{\Sigma}\}>0.$$
\end{proposition}
\begin{proof} Trivial as $\lambda: \Sigma\ni (S,K)\longmapsto \lambda(S,K)\in [0,1]$\ is continuous and $\widetilde{\Sigma}$\ is compact.

\end{proof}

\begin{corollary}\label{col: przeciecia transwersalne 2} Let $s,k,p,n\in \N$, $p<\min\{k,s\}$\ and let $\widetilde{\Sigma}$\ be a compact subset of
$$\Sigma=\{(S,K)\in\mathbb{G}_{s,n}\times\mathbb{G}_{k,n}:\quad \dim (S\cap K)=p\}.$$

Then there exists $C>0$\ such that for any linear subspaces $R,L$\ of $\R^{n}$\ and for any $(S,K)\in\widetilde{\Sigma}$
$$d(R\cap L, S\cap K)\leq C\cdot (d(R,S)+d(L,K)).$$
\end{corollary}

\begin{proof}
By Proposition \ref{prop: inf lambda} and \ref{prop: lambda pierwsza nierownosc} $ii)$\ the above inequality holds true with the constant
$$C=\frac{1}{\inf\{\lambda(S,K):\quad (S,K)\in\widetilde{\Sigma}\}}.$$
\end{proof}

\begin{definition}\label{def: Verdier condition}
Let $\Lambda$, $\Gamma$\ be $C^{2}$\ submanifolds of $\R^{n}$, $\Gamma\subset\overline{\Lambda}\setminus\Lambda$. We say that
the pair $(\Lambda,\Gamma)$\ satisfies \textit{the Verdier condition at }$x_{0}\in\Gamma$\ (notation: $\W^{V}(\Lambda, \Gamma, x_{0})$)\ if
there exists an open
 neighbourhood $U_{x_{0}}$\ of $x_{0}$\ in $\R^{n}$\ and $C_{x_{0}}>0$\ such that
 $$\forall x\in\Gamma\cap U_{x_{0}}\text{   }\forall y\in \Lambda\cap U_{x_{0}}\quad d\left(T_{x}\Gamma, T_{y}\Lambda\right)\leq C_{x_{0}}\cdot|x-y|.$$
 \vskip 1 mm
We say that $(\Lambda,\Gamma)$\ satisfies the Verdier condition (notation $\W^{V}(\Lambda,\Gamma)$)\ if for each $x_{0}\in\Gamma$\ we have $\W^{V}(\Lambda,\Gamma,x_{0})$.
\end{definition}

In 1998 Ta Le Loi proved that

\begin{theorem}
The Verdier condition is definable and generic.
\end{theorem}
\begin{proof} See \cite{TL1} (compare to \cite{LSW} and \cite{DW}).
\end{proof}

\begin{remark} As it was explained in \cite{Ver}, the Verdier condition is invariant under definable
$C^{q}$\ diffeomorphisms, $q\geq 2$. However, it is not $C^{1}$\ invariant, as was shown in \cite{BT}.
\end{remark}

Now we prove that the Verdier condition has the projection property with respect to weakly Lipschitz mappings of class $C^{q}$, where  $q\geq 2$.

\begin{theorem}
Let $q\geq 2$\ and let $\Lambda$, $\Gamma$\ be $C^{q}$\ submanifolds of $\R^{n}\times \R^{m}$, $\Gamma\subset\overline{\Lambda}\setminus\Lambda$\ and $\dim \Gamma<\dim \Lambda$\
 \footnote{Again we have to assume that $\{\Lambda,\Gamma\}$\ is a $C^{q}$\ stratification of $\Lambda\cup\Gamma$.}.
Let the mapping $\pi_{1}|_{\Lambda\cup\Gamma}$\ be a homeomorphic embedding such that
$\pi_{1}|_{\Lambda}$, $\pi_{1}|_{\Gamma}$\ are $C^{q}$\ embeddings\footnote{As before $\pi_{1}: \R^{n}\times\R^{m}\longrightarrow\R^{n}$\ and
$\pi_{2}: \R^{n}\times\R^{m}\longrightarrow\R^{m}$\ denote natural projections.}.
Assume that the mapping $\pi_{2}|_{\Lambda\cup\Gamma}\circ\left(\pi_{1}|_{\Lambda\cup \Gamma}\right)^{-1}$\
is weakly Lipschitz on the $C^{q}$\ stratification $\{\pi_{1}(\Lambda), \pi_{1}(\Gamma)\}$.

Then
$$\W^{V}(\Lambda,\Gamma)\Longrightarrow \W^{V}(\pi_{1}(\Lambda),\pi_{1}(\Gamma)).$$
\end{theorem}

\begin{proof} Put $x'=\pi_{2}|_{\Lambda\cup\Gamma}\circ (\pi_{1}|_{\Lambda\cup\Gamma})^{-1}(x)$, for any $x\in\pi_{1}(\Lambda\cup\Gamma)$.
\vskip 2 mm

Let $x_{0}\in \pi_{1}(\Gamma)$. After making a suitable $C^{2}$\ change of coordinates in a neighbourhood of $(x_{0}, x'_{0}$) we can assume that
$\pi_{1}(\Gamma)=\R^{k}\times \{0\}^{n-k}$\ and $\Gamma=\pi_{1}(\Gamma)\times \{0\}^{m}=\R^{k}\times \{0\}^{n+m-k}=\R^{k}$.
Then $x'=0$, for any $x\in\pi_{1}(\Gamma)$.

There exists an open neighbourhood $U_{(x_{0},0)}$\
of the point $(x_{0},0)$\ in $\R^{n}\times \R^{m}$ \ and a constant $C_{(x_{0},0)}>0$\ such that
$$d\left(T_{(x,0)}\Gamma, T_{(y,y')}\Lambda\right) \leq C_{(x_{0},0)}\cdot |(x,0)-(y,y')|,$$
for each $(x,0)\in \Gamma\cap U_{(x_{0},0)}$\ and $(y,y')\in \Lambda\cap U_{(x_{0},0)}$.

Since the mapping $\pi_{2}|_{\Lambda\cup\Gamma}\circ(\pi_{1}|_{\Lambda\cup\Gamma})^{-1}$\ is weakly Lipschitz on the stratification
$\{\pi_{1}(\Lambda), \pi_{1}(\Gamma)\}$,\ there is a neighbourhood $U_{x_{0}}$\ of the point $x_{0}$\ in $\R^{n}$\ and
a constant $L_{x_{0}}>0$\ such that
$$\frac{|0-y'|}{|x-y|}\leq L_{x_{0}},$$
for each $x\in \pi_{1}(\Gamma)\cap U_{x_{0}}$\ and $y\in \pi_{1}(\Lambda)\cap U_{x_{0}}$.

Without loss of generality we may assume that $U_{x_{0}}=U_{(x_{0},0)}\cap \left(\R^{n}\times\{0\}^{m}\right)$.
Then from the above argument, for all points $(x,0)\in \Gamma\cap U_{(x_{0},0)}$\ and $(y,y')\in \Lambda\cap U_{(x_{0},0)}$
$$d(T_{(x,0)}\Gamma, T_{(y,y')}\Lambda)
\leq C_{(x_{0},0)}\cdot |(x,0)-(y,y')|\leq C_{(x_{0},0)}\cdot \sqrt{1+(L_{x_{0}})^{2}}\cdot |x-y|,$$
in other words
$$d(\R^{k}, T_{(y,y')}\Lambda)\leq C_{(x_{0},0)}\cdot\sqrt{1+L^{2}_{x_{0}}}|x-y|.$$
Hence, after diminishing $U_{(x_{0},0)}$,
$$d(\R^{k},T_{(y,y')}\Lambda)\leq 1-\alpha,$$
where $0<\alpha<1$\ and then denoting by $K_{y}$\ the orthogonal projection of $\R^{k}$\ onto $T_{(y,y')}\Lambda$, we get a $k$\nolinebreak\ dimensional
subspace $K_{y}$\ of $T_{(y,y')}\Lambda$\ such that
$$d(\R^{k}, K_{y})=d(\R^{k}, T_{(y,y')}\Lambda).$$
Consequently, by Lemma \ref{lem: Lipsch_rzut_funkcji_d} $ii)$:
$$d\left(T_{x}\pi_{1}(\Gamma), T_{y}\pi_{1}(\Lambda)\right)=d\left(\R^{k}, \pi_{1}\left(T_{(y,y')}\Lambda\right)\right)=$$
$$d\left(\pi_{1}(\R^{k}\times\{0\}), \pi_{1}\left(K_{y}\right)\right)\leq C_{\alpha}\cdot d\left(\R^{k}, T_{(y,y')}\Lambda\right)\leq
\widetilde{C}\cdot |x-y|,$$
which completes the proof.

\end{proof}

\begin{corollary} The Verdier condition has the projection property with respect to weakly Lipschitz mappings of class $C^{q}$, where $q\geq 2$.
\end{corollary}

It remains to deal with the property of lifting the Verdier condition with respect to locally Lipschitz mappings of class $C^{q}$, $q\geq 2$. The
argument is analogous to that in the case of the Whitney (B) condition.
\vskip 3 mm

\begin{theorem}\label{trm: transversality Verdier} Theorem \ref{trm: Whitney transveral theorem} remains true when the Whitney (B) condition is replaced
by the Verdier condition, assuming that $q\geq 2$.
\end{theorem}
\begin{proof} Let $x_{0}\in\Gamma_{1}\cap\Gamma_{2}$. There exists a neighbourhood $U_{x_{0}}$\ of $x_{0}$\ and a constant $C_{x_{0}}>0$\ such that
for each $x\in\Gamma_{i}\cap U_{x_{0}}$\ and each $y\in\Lambda_{i}\cap U_{x_{0}}$
$$d\left(T_{x}\Gamma_{i}, T_{y}\Lambda_{i}\right)\leq C_{x_{0}}\cdot |x-y|,$$
where $i=1,2$.
\vskip 2 mm

\textit{Case I.} Assume that $\dim \Lambda_{1}\cap\Lambda_{2}=\min\{\dim \Lambda_{1},\dim\Lambda_{2}\}=\dim\Lambda_{1}$\ (the case, when
 $\dim\Lambda_{1}\cap\dim\Lambda_{2}=\min\{\dim \Lambda_{1},\dim\Lambda_{2}\}=\dim\Lambda_{2}$\ is similar).
 Then $\Lambda_{1}\cap \Lambda_{2}$\ is open in $\Lambda_{1}$, hence $T_{y}(\Lambda_{1}\cap\Lambda_{2})=T_{y}\Lambda_{1}$\ and
 $$d\left(T_{x}(\Gamma_{1}\cap\Gamma_{2}), T_{y}(\Lambda_{1}\cap \Lambda_{2})\right)\leq d\left(T_{x}\Gamma_{1},T_{y}\Lambda_{1}\right)\leq C_{x_{0}}\cdot |x-y|$$
 for each $x\in\Gamma_{1}\cap\Gamma_{2}\cap U_{x_{0}}$\ and $y\in\Lambda_{1}\cap\Lambda_{2}\cap U_{x_{0}}$,
which completes the proof of the Case.
\vskip 2 mm
\textit{Case II.} Assume now that $\dim \Lambda_{1}\cap\Lambda_{2}<\min\{\dim \Lambda_{1},\dim\Lambda_{2}\}$. Let
$$\Sigma=\{(S, K)\in\mathbb{G}_{\dim\Lambda_{1},n}\times\mathbb{G}_{\dim\Lambda_{2},n}: dim (S\cap K)=\dim\Lambda_{1}\cap\Lambda_{2}\}.$$
Then by the assumptions, after perhaps diminishing the neighbourhood $U_{x_{0}}$\ the closure $\widetilde{\Sigma}$\ of the set
$$\{(T_{y}\Lambda_{1}, T_{y}\Lambda_{2})\in\mathbb{G}_{\dim \Lambda_{1},n}\times \mathbb{G}_{\dim \Lambda_{2},n}:\quad y\in\Lambda_{1}\cap\Lambda_{2}\cap U_{x_{0}}\}$$
is a closed subset of $\Sigma$.
By Corollary \ref{col: przeciecia transwersalne 2} there exists a constant $C>0$\ such that for all $x\in\Gamma_{1}\cap \Gamma_{2}\cap U_{x_{0}}$,
$y\in\Lambda_{1}\cap \Lambda_{2}\cap U_{x_{0}}$\ we have
$$d\left(T_{x}\Gamma_{1}\cap T_{x}\Gamma_{2},T_{y}\Lambda_{1}\cap T_{y}\Lambda_{2}\right)\leq C\cdot\left(d\left(T_{x}\Gamma_{1},T_{y}\Lambda_{1}\right)+d\left(T_{x}\Gamma_{2},T_{y}\Lambda_{2}\right)\right).$$
Consequently,
$$d\left(T_{x}\left(\Gamma_{1}\cap \Gamma_{2}\right),T_{y}\left(\Lambda_{1}\cap \Lambda_{2}\right)\right)\leq C\cdot 2 C_{x_{0}}\cdot |x-y|.$$
\end{proof}

\begin{theorem} Theorem \ref{trm: Whitney (B) lifting property} remains true when the Whitney (B) condition is replaced by the Verdier condition, assuming
that $q\geq 2$.
\end{theorem}
\begin{proof} It follows from Theorem \ref{trm: transversality Verdier} in the same way as Theorem \ref{trm: Whitney (B) lifting property} follows from
Theorem \ref{trm: Whitney transveral theorem}.

\end{proof}

\begin{corollary} The Verdier condition is a $\WL$\ condition of class $C^{q}$, $q\geq 2$.
\end{corollary}

\begin{corollary} Theorem \ref{trm: Q invariance} holds true for the Verdier condition with $q\geq 2$.
\end{corollary}

\textbf{Acknowledgements.} I would like to express my profound gratitude to Wies\l aw Paw\l ucki for his great tutorial, friendly support and
many valuable discussions.
%%%%%%%%%%%%%%%%%%%%%%%%%%%%%%%%%%%%%%%%%%%%%%%%%%
\medskip
%%%%%%%%%%%%%%%%%%%%%%%%%%%%%%%%%%%%%%%%%%%%%%%%%%
%References
%%%%%%%%%%%%%%%%%%%%%%%%%%%%%%%%%%%%%%%%%%%%%%%%%%
%\bibliographystyle{amsplain}

\end{document}